\documentclass{amsart}
\usepackage{amsmath,amssymb,amsthm,amscd, mathtools, latexsym}
\usepackage{enumerate,varioref, dsfont, fancyhdr}
\usepackage{tikz-cd}
\usepackage{enumitem}
\usepackage{multicol}
\usepackage{hyperref}
\usepackage{url}
\usepackage{textcomp}

\newtheorem{Thm}{Theorem}[section]
\newtheorem{Lem}{Lemma}[section]
\newtheorem{Prop}{Proposition}[section]

\newtheorem{Cor}{Corollary}[section]

\newtheorem{Rem}{Remark}[section]

{\theoremstyle{plain}

\newtheorem*{Ass*}{Assumption}

}
\renewcommand{\ker}{\text{ker}}

\newcommand{\Sym}{\text{Sym}}
\newcommand\restr[2]{\ensuremath{\left.#1\right|_{#2}}}

\def\cit{{\mathbb C}}

\def\qit{{\mathbb Q}}

\def\pit{{\mathbb P}}

\def\0{{\mathcal O}}

\def\End{\mathop{\rm End}\nolimits}

\def\h{{\mathfrak h}}

\def\F{{\mathcal F}}
\def\h{{\mathfrak h}}

\def\M{{\mathcal M}}
\def\N{{\mathcal N}}

\def\E{{\mathcal E}}

\def\X{{\mathcal X}}

\begin{document}
\title[The Chow ring of a cubic hypersurface]{The Chow ring of a cubic hypersurface}

\author{H. Anthony Diaz}

\begin{abstract} We study the product structure on the Chow ring (with rational coefficients) of a cubic hypersurface in projective space and prove that the image of the product map is as small as possible.

\end{abstract}

\maketitle

\section{Introduction}

\noindent Let $X$ be a smooth projective variety over $\cit$ and let $A^{*} (X)$ denote the Chow ring of $X$ with rational coefficients, graded by codimension. Describing the product structure of $A^{*} (X)$ is a difficult problem. When the total cycle class map is injective (e.g., when $X$ is homogeneous), the product structure on $A^{*}(X)$ is determined by that of the cohomology ring of $X$. On the other hand, it is sometimes possible to describe the product structure in certain special cases where the cycle class map is no longer injective. For instance, there are the two famous examples of Abelian varieties and $K3$ surfaces. When $X$ is an Abelian variety, Beauville \cite{B1} showed that there is a natural multiplicative decomposition:
\[ A^{i} (X) = \bigoplus_{j} A^{i}_{(j)} (X)\]
for which it is expected that $A^{i}_{(j)} (X)= 0$ for $j<0$ and that the cycle class map
\[ A^{i}_{(0)} (X) \to H^{2i} (X, \qit)\]
is injective (thus giving a splitting to the conjectural Bloch-Beilinson conjecture; see, for instance, \cite{B}). Another well-known case is when $X$ is a complex projective $K3$ surface
which was studied by Beauville and Voisin \cite{BV}; in particular, they showed that the image of the map
\[ A^{1} (X) \otimes A^{1} (X) \xrightarrow{\cdot} A^{2} (X) \]
has rank $1$, a surprising discovery considering that $A^{2} (X)$ is not representable by Mumford's theorem \cite{M1}.\\
\indent Beyond these fundamental cases, it is natural to wonder to what extent one may determine the structure of the Chow ring of (for instance) a hypersurface of degree at least $3$. For small degree hypersurfaces, it is expected that the Chow group will be of rank $1$ when the codimension is either very small or very large (see, for instance, \cite{Mb}). On the other hand, since hypersurfaces of degree at least $3$ possess transcendental cohomology in the middle degree, a generalization of Mumford's theorem (\cite{V4} Theorem 3.20) shows that the cycle class map is not injective in every degree. Thus, the total Chow group is a mysterious invariant, even for a cubic hypersurface. In this note, our aim will be to prove the following explicit characterization of the product structure for the Chow ring of a cubic hypersurface: 
\begin{Thm}\label{int-prod} Let $X \subset \pit^{n+1}_{\cit}$ be a smooth hypersurface of degree $3$. Then, the image of the product map
\[ A^{i} (X) \otimes A^{j} (X) \xrightarrow{\cdot} A^{i+j} (X) \]
has rank $1$ for $i,j>0$.\footnote{Note this is only interesting in the case $i+j<n$, since $X$ is rationally connected.}
\end{Thm}
\noindent We give an overview of the method of proof for Theorem \ref{int-prod}. In the spirit of \cite{V3} and \cite{Fu} we give a decomposition of the small diagonal on $X^{3}$. Since $A^{*} (X^{3})$ is potentially vast, we work with the a priori finite rank subspace $R^{*} (X^{3})$ genered by $A^{1} (X^{3})$ and the pull-backs of the diagonal classes. Since the small diagonal admits a decomposition cohomologically (by the K\"unneth theorem) for any hypersurface, it will suffice to prove that the cycle class map
\begin{equation} R^{*} (X^{3}) \hookrightarrow H^{*} (X^{3}, \qit)\label{hook-cyc} \end{equation}
We should note that a similarly defined ring $R^{*}$ was considered by Yin in \cite{Y} for (arbitrarily many) products of $K3$ surfaces. He showed that, when $X$ is a $K3$ surface, the injectivity of the cycle class map restricted to $R^{*}(X^{n})$ is equivalent to the finite-dimensionality of the Chow motive of $X$ \cite{K}. On the other hand, for hypersurfaces of general type, one cannot expect results analogous to (\ref{hook-cyc}) (see, for instance, \cite{OG}).\\ 
\indent As in \cite{V3} and \cite{Fu}, the proof will also involve the Fano variety of lines, which is better understood for cubic hypersurfaces. A thorough study of the Fano variety of lines $F$ on an $n$-dimensional cubic hypersurface was completed by Altman and Kleiman in \cite{AK}, in which they show (among other things) that $F$ is smooth and has dimension $2(n-2)$. In the case of $n=3$, $F$ is a surface of general type which possesses a great many remarkable properties, used by Clemens and Griffiths (in their well-known paper \cite{CG}) to establish that the smooth cubic threefold is not rational. For $n=4$, Beauville and Donagi showed in \cite{BD} that $F$ has the structure of a hyper-K\"ahler variety that is deformation equivalent to the second punctual Hilbert scheme of a $K3$ surface. For $n \geq 5$, $F$ is rationally connected. Quite recently, Galkin and Schinder in \cite{GS} were able to prove the following relation in the Grothendieck ring of varieties:
\[ [X^{[2]}] = [X][\pit^{n+1}] + \mathds{L}^{2}[F] \in K_{0} (Var) \]
where $X^{[2]}$ denotes the second punctual Hilbert scheme of $X$. A motivic version of their result was then obtained by Laterveer in \cite{L}; we will use the motivic version of this relation to reduce the problem of proving (\ref{hook-cyc}) to proving an analogous result for $F$ and eventually $F \times X$, which it turns out are easier to prove. We should note that with the exception of the Galkin-Shinder result, the techniques used here are quite elementary and do not rely on the fact that the degree of $X$ is $3$. It is likely that one could produce similar (albeit weaker) results for other rationally connected hypersurfaces in projective space, providing that one can establish an analogous motivic relation with the Chow motive of the corresponding Fano variety of lines.
\subsection*{Acknowledgements} 
The author would like to thank Lie Fu and Charles Vial for their interest in this result.
\section{Preliminaries}
\subsection*{Conventions and Notation}
Unless otherwise specified, we will take our ground field to be $\cit$. $A^{*}$ will denote the Chow ring with rational coefficients and $H^{*} (-,\qit)$ will denote singular cohomology with $\qit$ coefficients. We will also fix the following notation from now on:
\begin{itemize}
\item $X$, a smooth cubic hypersurface in $\pit^{n+1}$;
\item $X^{[2]}$, the second punctual Hilbert scheme on $X$;
\item $E_{X}$, the exceptional divisor on $X^{[2]}$;
\item $h:=c_{1}(\mathcal{O}_{X}(1))$;
\item $F := F(X)$, the Fano variety of lines on $X$;
\item $G_{n+1}:= Gr(\pit^{1}, \pit^{n+1})$, the Grassmanian of lines in $\pit^{n+1}$;
\item $c_{i}:= c_{i}(\E) \in CH^{i} (G_{n+1})$, the $i^{th}$ Chern class of $\E$ (when there is no ambiguity, we will also denote by $c_{i}$ the corresponding Chern class $\restr{c_{i}}{F}$ on $F$); 
\item $\iota: F \hookrightarrow G_{n+1}$, the canonical imbedding;
\item $\mathcal{E}$, the tautological rank $2$ vector bundle on $G_{n+1}$ (by abuse, we also let $\E$ denote the restriction $\iota^{*}\E$;
\item $j: P:=\pit(\E) \hookrightarrow F \times X$, the universal line over $F$;  
\item $p_{F}:  F \times X \to F$ and $p_{X}: F \times X\to F$, the corresponding projections;
\item $p:= \restr{p_{F}}{P}$, $q:= \restr{p_{X}}{P}$.
\end{itemize}
\subsection{Tautological rings}
We will  consider the following {\em tautological rings}:
\begin{itemize}
\item $R^{*}(X) \subset A^{*} (X)$, the subring generated by $A^{1} (X)$;
\item $R^{*}(F) \subset A^{*} (F)$, the subring generated by $c_{1}$ and $c_{2}$;
\item $R^{*} (X \times X) \subset A^{*} (X\times X)$, the subring generated by $A^{1} (X \times X)$ and $\Delta_{X}$ (or, equivalently, by $\pi_{1}^{*}R^{*}(X)$, $\pi_{2}^{*}R^{*}(X)$ and $\Delta_{X}$);
\item $R^{*}(F \times X) \subset A^{*} (F \times X)$, the subring generated by $p_{F}^{*}R(F)$, $p_{X}^{*}R(X)$ and $\Gamma$.
\end{itemize}
\noindent Moreover, if $Y$ is any of the above varieties and $P_{Y} \xrightarrow{\pi} Y$ is a projective bundle, we take the tautological ring $R^{*}(P_{Y})$ to be the subring of $A^{*}(P_{Y})$ generated by $\pi^{*}A^{*}(Y)$ and $c_{1} (\mathcal{O}_{P_{Y}}(1))$.\\
\indent We would like to give a characterization of all of the above tautological rings. To this end, we introduce some more set-up. Indeed, let 
$\pit=\pit(H^{0} (\pit^{n+1}, \mathcal{O}_{\pit^{n+1}} (3))$ be the projective space parametrizing cubic hypersurfaces in $\pit^{n+1}$ and let $\X \to \pit$ be the universal family of cubics over $\pit$. Let $\mathcal{F} \to \pit$ be the corresponding family of Fano varieties of lines. Furthermore, let $\eta$ be the generic point of $\pit$ and let $\X_{\eta}$ and $\F_{\eta}$ be the corresponding generic fibers. It is easy to see (for instance) that
\begin{equation}\label{incl}
\begin{split}
R^{k}(X) & \subset \text{Im}\{A^{k} (\X_{\eta}) \to A^{k} (X) \}\\
R^{k} (F) & \subset \text{Im}\{A^{k} (\F_{\eta}) \to A^{k} (\F_{\eta})\}\\
R^{k} (X \times X) & \subset \text{Im}\{A^{k} (\X_{\eta} \times_{\eta} \X_{\eta}) \to A^{k} (X \times X)\}
\end{split}
\end{equation}
\begin{Lem}\label{first} The inclusions in (\ref{incl}) are all equalities.
\begin{proof} We first prove the first two inclusions. To this end, note that there are natural inclusions:
\[\begin{split}
\mathcal{X} & \subset \pit \times \pit^{n+1}\\
\mathcal{F} & \subset \pit \times G_{n+1} 
\end{split}\]
and let $\pi_{\pit^{n+1}}: \mathcal{F} \to \pit^{n+1}$ and $\pi_{G_{n+1}}: \mathcal{F} \to G_{n+1}$ be the projections. Then, we observe that both $\pi_{\pit^{n+1}}$ and $\pi_{G_{n+1}}$ are projective bundles. Indeed, given any $p \in \pit^{n+1}$ (resp., $\ell \in  G_{n+1}$), the set of cubics in $\pit$ containing $p$ (resp., $\ell$) is a projective linear subspace, and the dimension of this subspace does not depend on $p$ (resp., $\ell$), from which it follows that these are projective bundles. From the projective bundle formula \cite{F} Theorem 3.3, we have
\[\begin{split}
A^{k} (\X) & = \bigoplus_{j\leq k} \pi^{*}A^{j} (\pit^{n+1})\cdot H^{k-j}\\
A^{k} (\F) & = \bigoplus_{j\leq k} \pi^{*}A^{j} (G_{n+1})\cdot H^{k-j}
\end{split}\]
where $H = c_{1} (\mathcal{O}_{\pit}(1))$. We would like to show that
\[\begin{split}
R^{k} (X) &=  \text{Im}\{A^{k} (\X) \to A^{k} (X) \}\\
R^{k} (F) &=  \text{Im}\{A^{k} (\F) \to A^{k} (F) \}
\end{split}\]
for which it suffices to show that the subspaces 
\begin{equation}\pi^{*}A^{j} (\pit^{n+1})\cdot H^{k-j}, \ \pi^{*}A^{j} (G_{n+1})\cdot H^{k-j}\label{vanish-sub} \end{equation} 
vanish upon restriction to a fiber vanishes for $j<k$. Letting $\pit' \subset \pit$ be a divisor whose class in $A^{1} (\pit)$ is $H$. By the projection formula, it follows that the subspaces in (\ref{vanish-sub}) are supported over $\pit'$ for $j<k$. The desired vanishing then follows.\\
\indent To prove that the final inclusion in (\ref{incl}) is an equality, we let $S \subset \pit$ denote the subspace of cubic hypersurfaces which are smooth and we would like to show that
\begin{equation} R^{k} (X \times X) =  \text{Im}\{A^{k} (\X \times_{S} \X) \to A^{k} (X \times X) \}\label{third}\end{equation}
In this case, there is a natural inclusion:
\[ \X \times_{S} \X \subset \pit \times \pit^{n+1} \times \pit^{n+1} \]
Set $Y = \pit^{n+1} \times \pit^{n+1}$ and as in \cite{V} Lemma 3.13 one can consider the corresponding inclusion of blow-ups along the diagonal 
\[ \widetilde{X \times_{S} \X} \subset \pit \times \tilde{Y}\]
and it follows from the proof of loc. cit. that the corresponding projection $\widetilde{X \times_{S} \X}$ is an open subset of a projective bundle $P_{\tilde{Y}}$ over $\tilde{Y}$. As in the previous paragraph, one shows that the pull-back
\begin{equation} A^{k} (\tilde{Y}) \to A^{k}(\widetilde{X \times_{S} \X})\label{pulled} \end{equation}
is surjective. Indeed, the tautological line bundle on $P_{\tilde{Y}} \subset \pit \times \tilde{Y}$ is pulled back from $\pit$ and hence vanishes on $\widetilde{X \times_{S} \X}$ since $A^{1} (S) =0$. The same argument as in the previous paragraph then applies to show that all the summands of $A^{*} (P_{\tilde{Y}})$ except $A^{*}(\tilde{Y})$ vanish on $A^{*}(\widetilde{X \times_{S} \X})$. The surjectivity of (\ref{pulled}) then follows. From \cite{F} Prop. 6.7, $A^{*} (\tilde{Y})$ is generated as a ring by
$\epsilon^{*}A^{*} (Y)$ and the image of the push-forward $A^{*-1} (E_{Y}) \to A^{*} (\tilde{Y})$, where $\epsilon: \tilde{Y} \to Y$ is the blow-up and $E_{Y}$ is the exceptional divisor. Since 
\[ A^{*} (Y) = A^{*} (\pit^{n+1}) \otimes  A^{*} (\pit^{n+1})\]
it follows that $A^{*} (\widetilde{X \times_{S} \X})$ is generated as a ring by $\epsilon_{i}^{*}A^{*} (\X)$ (where $\epsilon_{i}: \widetilde{X \times_{S} \X} \to \X \times_{S} \X \xrightarrow{\pi_{i}} \X$ are the natural projections) and the image of the push-forward
\[A^{*-1} (E_{\X \times_{S} \X}) \to A^{*} (\widetilde{X \times_{S} \X}) \]
By the projection formula, it follows that
\[ A^{*} (\X \times_{S} \X) = \epsilon_{*}A^{*} (\widetilde{X \times_{S} \X})\]
is generated as a ring by $\pi_{i}^{*}A^{*} (\X)$ and the image of the push-forward
\[A^{*-n} (\X) \xrightarrow{\Delta_{\X*}} A^{*} (\widetilde{X \times_{S} \X}) \]
(Recall that $E_{\X \times_{S} \X}$ lies over $\Delta_{\X}$.) But since $\Delta_{X}^{*}$ is surjective, the image of $\Delta_{X*}$ coincides with 
\[A^{*} (\widetilde{X \times_{S} \X})\cdot \Delta_{\X} \]
This gives (\ref{third}), as desired.
\end{proof}
\end{Lem}
\noindent In the case of $R(X)$ and $R(X \times X)$, there is also the following basic result, which we will eventually prove for the other tautological rings.
\begin{Lem}\label{basic} The cycle class maps
\[ R^{2k}(X) \to H^{2k} (X, \qit(k)), \ R^{2k}(X \times X) \to H^{2k} (X\times X, \qit(k))\]
are injective for all $k$.
\begin{proof} The injectivity of the first is trivial, as is the injectivity of the second for $k<n$. In general, we note that an additive generating set for $R^{*} (X \times X)$ is given by:
\begin{equation} \pi_{1}^{*}h^{r_{1}}\cdot\pi_{2}^{*}h^{r_{2}}, \  \Delta_{X}\cdot \pi_{1}^{*}h^{s_{1}}\cdot\pi_{2}^{*}h^{s_{2}}, \Delta_{X}^{2}\label{gen-set-2} \end{equation}
where $0 \leq r_{i}, s_{i} \leq n$. Now, it is easy to show using the K\"unneth theorem for cohomology that modulo homological($=$ numerical) equivalence 
\[\Delta\cdot \pi_{1}^{*}h^{s_{1}}\cdot\pi_{2}^{*}h^{s_{2}}\]
is a linear combination of cycles of the form $\pi_{1}^{*}h^{r_{1}}\cdot\pi_{2}^{*}h^{r_{2}}$ whenever $s_{1}, s_{2} \neq 0$. However, this also holds on the level of rational equivalence; indeed, if $s_{1} \geq 1$, we let $i: X \hookrightarrow \pit^{n+1}$ be the inclusion and compute
\begin{equation}\begin{split} \Delta_{X}\cdot \pi_{1}^{*}h^{s_{1}} &= \Delta_{X*}h^{s_{1}}\\
&= \frac{1}{3}\cdot\Delta_{X*}(i^{*}c_{1}(\mathcal{O}_{\pit^{n+1}}(1))^{s_{1}-1}\cdot 3h)\\
&= \frac{1}{3}\cdot(i \times i)^{*}\Delta_{\pit^{n+1}*}c_{1}(\mathcal{O}_{\pit^{n+1}}(1))^{s_{1}-1}\\
\end{split}\label{obvious}
\end{equation}
where the final line follows from the excess intersection formula \cite{F} \S 6.3. It is clear that the final line of (\ref{obvious}) is a linear combination of cycles of the form $\pi_{1}^{*}h^{r_{1}}\cdot\pi_{2}^{*}h^{r_{2}}$. The analogous statement is also true for $\Delta_{X}\cdot \pi_{2}^{*}h^{s_{2}}$ when $s_{2} \geq 1$. Finally, we observe that $\Delta_{X}^{2}$ is rationally equivalent to a multiple of $\pi_{1}^{*}h^{n}\cdot\pi_{2}^{*}h^{n}$, since $X$ is rationally connected. Thus, the generating set (\ref{gen-set-2}) truncates to a basis
\begin{equation}\pi_{1}^{*}h^{r}\cdot\pi_{2}^{*}h^{s}, \ \Delta_{X}\label{basis-2} \end{equation}
(with $0 \leq r, s \leq n$) for $R^{*} (X \times X)$. This is also a basis for the image of the above cycle class map; hence, the lemma.
\end{proof}
\end{Lem}

\subsection{Galkin-Shinder relation}
\noindent We would like to derive some relationships among the various tautological rings. To this end, we will need a result of Galkin and Shinder \cite{GS}, which relates a cubic hypersurface to its Fano variety of lines. To introduce it, we observe that there is a rational map:
\[ X^{[2]} \dashrightarrow \pit(\restr{T\pit^{n+1}}{X})\]
that sends a general $x+x'$ to the third point $x''$ on the line which passes through $x$ and $x'$ (along the exceptional divisor $E_{X}$, this line is the tangent line through $x$ in a given direction). This map is in fact birational and induces an isomorphisms between the open subsets $U \subset X^{[2]}$ and $U' \subset \pit(\restr{T\pit^{n+1}}{X})$, whose respective complements are given by 
\[ 
\begin{split}
Z & = \{ x+y \in X^{[2]} \ | \ \exists \ell \in F \text{ such that } x, y \in \ell \}\\
Z' & = \{(z \in \ell) \in \pit(\restr{T\pit^{n+1}}{X}) \ | \ \ell \subset X \}
\end{split}
\]
It is not difficult to see that $Z$ and $Z'$ are projective bundles over $F$ and, in fact, that 
$Z \cong \pit(\Sym^{2}(\E))$ and $Z' \cong \pit(\E)$. 
\noindent Using these observations, the authors of \cite{GS} obtain the following relation in the Grothendieck ring of varieties over $k$, $K_{0} (Var/k)$:
\[[X^{[2]}] = [X][\pit^{n+1}] + \mathds{L}^{2}[F]  \]
where $\mathds{L}$ denotes the class of $\mathit{A}^{1}$. By resolving the indeterminacy of this map, Laterveer obtains the following motivic version of this:
\begin{Thm}[\cite{L} Theorem 5]\label{Later} With the notation above, there is an isomorphism in $\M_{k}$:
\[ \h(F)(-2) \oplus \bigoplus_{i=0}^{n} \h(X)(-i) \cong \h(X^{[2]}) \]
In particular, $\h(F)$ is a summand of $\h(X^{[2]})(2)$.
\end{Thm}
\begin{Cor}\label{cyc-inj} The cycle class map
\[ R^{k} (F) \to H^{2k} (F, \qit(k) \]
is injective for all $k$.
\begin{proof} By Lemma \ref{first}, it suffices to show that rational equivalence and numerical equivalence agree for cycles on $\F_{\eta}$. (A specialization argument then gives the statement of the corollary.) By the previous theorem, we are reduced to proving the corresponding statement  for $\X_{\eta}$ and $\X_{\eta}^{[2]}$, but this follows from Lemma \ref{basic}.
\end{proof}
\end{Cor}
\begin{Cor}\label{non-deg} The intersection product
\[ R^{k} (F) \otimes R^{2(n-2)-k} (F) \to \qit \]
is non-degenerate.
\begin{proof} This follows from the proof of Corollary \ref{cyc-inj}, which in fact shows that rational equivalence and numerical equivalence coincide for cycles in $R^{*}(F)$.
\end{proof}
\end{Cor}
\begin{Rem} The proof of Corollary \ref{cyc-inj} shows that any numerically trivial polynomial in $c_{1}$ and $c_{2}$ is also rationally trivial. We should note that Theorem \ref{Later} is not really needed; all one essentially needs is the statement that the pull-back map
\[ A^{k} (\X_{\eta}^{[2]}) \to A^{k} (\mathcal{Z}_{\eta}) \]
is surjective, where $\mathcal{Z}_{\eta} \to \F_{\eta}$ is the $\pit^{2}$-bundle described above. This is rather straightforward, given the above characterization of $A^{*} (\F_{\eta})$.
\end{Rem}
\begin{Cor}\label{push-inj} The push-forward 
\[ \iota_{*}: R^{k} (F) \to A^{k+4} (G_{n+1})\]
is injective for all $k$.
\begin{proof} By definition $R^{k} (F) = \iota^{*}A^{k} (G_{n+1})$, so this follows from Corollary \ref{non-deg}.
\end{proof}
\end{Cor}
\noindent The following result gives an ``extra" relation among the Chern classes restricted to $F$; i.e., a relation which does not hold on the level of $A^{*}(G_{n+1})$.
\begin{Prop}\label{non-triv} There exists $P(x,y) \in \qit[x,y]$ weighted of degree $n-1$ for which
\[ P(c_{1}, c_{2}) = 0 \in R^{n-1}(F)\]
and for which the coefficient of $c_{1}^{n-1}$ is non-zero.
\begin{proof} The class $[F] \in A^{4} (G_{n+1})$ is of the form $Q(c_{1}, c_{2})$, where $Q(x,y)$ is a weighted degree $4$ polynomial; in fact, by \cite{F} Example 14.7.13 $Q(x,y)$ is a linear combination of $x^{2}y$ and $y^{2}$. Now, for dimension reasons the map
\[ A^{n-1} (G_{n+1}) \xrightarrow{\cdot [F]} A^{n+3} (G_{n+1}) \]
is not injective (note $G_{n+1}$ has dimension $2n$). It follows that there exists a non-zero weighted degree $n-1$ polynomial $P(x,y) \in \qit[x,y]$ for which 
\[ P(c_{1},c_{2})\cdot[F] = 0 \in A^{n+3} (G_{n+1}) \]
By Corollary \ref{push-inj} this means that
\[ \iota^{*}P(c_{1}, c_{2}) = 0 \in A^{n-1} (F) \]
We would like to show that the $x^{n-1}$ coefficient of $P(x,y)$ is non-zero. Assume by way of contradiction that instead
\[ P(x,y) = a_{0}x^{n-3}y + a_{1}x^{n-5}y^{2} + \ldots \in \qit[x,y]\]
for $a_{k} \in \qit$. Then, by the above characterization of $Q$, we have   
\begin{equation} R(x,y) = P(x,y)Q(x,y) = b_{0} x^{n-1}y^{2} +b_{1}x^{n-3}y^{3}+\ldots \in \qit[x,y]\label{bad-form} \end{equation}
for $b_{k} \in \qit$. Now, the intersection product structure on $A^{*} (G_{n+1})$ is well-known (see, for instance, \cite{F} p. 270). As a ring $A^{*} (G_{n+1})$ is the quotient of the ring of symmetric polynomials in $2$ variables modulo the ideal $I= (R_{n}(s_{1}, s_{2}), R_{n+1}(s_{1}, s_{2}))$, where $s_{1}$ and $s_{2}$ are the elementary symmetric polynomials in two variables and $R_{k}(x,y)$ is a weighted degree $k$ polynomial such that $R_{k} (s_{1}, s_{2})$ is the complete symmetric polynomial of degree $k$. Since $R(c_{1}, c_{2}) = 0$, $R(s_{1}, s_{2}) \in I$. This implies that there exist $m_{1}, m_{2}, m_{3}, m_{4} \in \qit$ such that
\[ R(x,y) = (m_{1}x^{2} + m_{2}y)R_{n+1}(x,y) + (m_{3}xy +m_{4}x^{3})R_{n} (x,y)\]
But given (\ref{bad-form}), this implies that $m_{1}=m_{2}=m_{3}=m_{4}=0$, which contradicts the fact that $P(x,y) \neq 0$. The result now follows.
\end{proof}
\end{Prop}

\subsection{Tautological ring of $F \times X$}

\noindent Now, we would like to give a characterization of $R^{*} (F \times X)$. To this end, we let $\mathcal{D}R^{*} (F \times X)$ denote the subring of {\em externally decomposable cycles}; i.e., the subring generated by $p_{F}^{*}R^{*}(F)$ and $p_{X}^{*}R^{*}(X)$. We begin with the following observation that follows directly from Corollary \ref{cyc-inj} and the K\"unneth theorem in cohomology:
\begin{Lem}\label{Dr-inj} The cycle class map
\[ \mathcal{D}R^{k} (F \times X) \to H^{2k} (F \times X, \qit(k))\]
is injective for all $k$.
\end{Lem}
\noindent Now, let $\mathcal{D}R^{*} (F \times X)'$ be the subspace of $R^{*} (F \times X)$ having as a basis:
\begin{equation} \{ \Gamma, \Gamma\cdot p_{F}^{*}c_{1}, \ldots \Gamma\cdot p_{F}^{*}c_{1}^{n-2} \}\label{base}\end{equation}
(that this is linearly independent follows because of degree reasons).
\begin{Thm}\label{key} $R^{*} (F \times X) = \mathcal{D}R^{*}(F \times X) \oplus \mathcal{D}R^{*} (F \times X)'$.
\begin{proof} We first show that
\[ R^{*} (F \times X) = \mathcal{D}R^{*}(F \times X) + \mathcal{D}R^{*} (F \times X)'\]
To this end, we note that the tautological ring $R^{*}(F \times X)$ is generated as a ring by $DR(F \times X)$ and $\Gamma$. Thus, it suffices to show that 
\[ \Gamma\cdot R(F \times X) \subset \mathcal{D}R^{*}(F \times X) + \mathcal{D}R^{*} (F \times X)' \]
This is accomplished in the three results below.
\begin{Lem}\label{kunn} $\Gamma \cdot p_{X}^{*}h, \Gamma\cdot p_{F}^{*}c_{2}, \Gamma \cdot p_{F}^{*}c_{1}^{n-1} \in \mathcal{D}R^{*} (F \times X)$. 
\begin{proof} We first prove this for $\Gamma \cdot p_{X}^{*}h$. By the proof of \cite{V2} Prop. 3.3 (see also \cite{SV} Prop. A.6), there is a decomposition:
\[ \Delta_{X} = \Delta_{0} + \Delta_{1} \in A^{n} (X \times X) \]
where 
\[ \Delta_{0} = \sum_{j} a_{j}\pi_{1}^{*}h^{j}\pi_{2}^{*}h^{n-j} \]
and where $\Delta_{0}\cdot\pi_{2}^{*}h = 0$. Now, we have
\[\begin{split}\Gamma\cdot p_{F}^{*}h & = (p_{F} \times \text{id}_{X})_{*}(p_{X} \times \text{id}_{X})^{*}(\Delta_{X}\cdot\pi_{2}^{*}h)\\
& = (p_{F} \times \text{id}_{X})_{*}(p_{X} \times \text{id}_{X})^{*}(\Delta_{0}\cdot\pi_{2}^{*}h)\\
& = \sum_{j} a_{j}\cdot p_{F}^{*}(\Gamma^{*}(h^{j}))\cdot p_{X}^{*}h^{n-j+1}  \end{split}\]
which lies in $\mathcal{D}R^{*} (F \times X)$.\\
\indent  For $\Gamma\cdot p_{F}^{*}c_{2}$, we observe that
\[ p_{X}^{*}h^{2} = p_{F}^{*}c_{1}\cdot p_{X}^{*}h - p_{F}^{*}c_{2} \in A^{2} (P) \]
Intersecting with $\Gamma$ then gives the desired result, by the previous verification. Finally, the statement for $\Gamma \cdot p_{F}^{*}c_{1}^{n-1}$ follows from Proposition \ref{non-triv} and the previous two verifications.
\end{proof}
\end{Lem}
\begin{Cor}\label{cor-kunn} $\Gamma\cdot \mathcal{D}R^{*} (F \times X) \subset \mathcal{D}R^{*} (F \times X) + \mathcal{D}R^{*} (F \times X)'$
\end{Cor}
\begin{Lem} $\Gamma^{2} \in  \mathcal{D}R^{*} (F \times X) + \mathcal{D}R^{*} (F \times X)'$
\begin{proof} We begin by noting that
\[ \Gamma^{2} = j_{*}(c_{n-1}(\N)) \in A^{2n-2} (F \times X) \]
where $\N$ is the vector bundle on $F \times X$ fitting into the adjunction exact sequence:
\[ 0 \to TP \to j^{*}(TF \oplus TX) \to \N \to 0 \]
We retain the notation of the previous proof and there are exact sequences:
\begin{equation} \begin{split} 0 \to p^{*}TF \to TP \to \mathcal{Q} \to 0\\
0 \to \mathcal{O}_{P} \to p^{*}\E \otimes \mathcal{O}_{P}(1) \to \mathcal{Q} \to 0
\end{split}\label{exact} \end{equation}
Since $c_{1}(\mathcal{O}_{P}(1)) = q^{*}h$, it follows from Cartan's formula and induction that $c_{i} (TP)$ is in the $\qit$-algebra generated by 
\[\{ p^{*}c_{1}, p^{*}c_{2}, q^{*}h \} \subset j^{*}\mathcal{D}R^{*} (F \times X)\] 
Again, using Cartan's formula and induction, it follows that $c_{i} (\N)$ is also. It follows that
\[ \Gamma^{2} \in \Gamma\cdot \mathcal{D}R^{*} (F \times X) \]
By Corollary \ref{cor-kunn}, this gives the desired result.
\end{proof}
\end{Lem}

\noindent What remains now is to show that $\mathcal{D}R^{*} (F \times X) \cap \mathcal{D}R^{*} (F \times X)' = 0$, which can be verified cohomologically. Indeed, we let
\[H^{n}_{prim} (X, \qit) := \ker \{ \cup h: H^{n} (X, \qit) \to H^{n+2} (X, \qit(1)) \} \]
Then, we have the following lemma:
\begin{Lem} \label{second} $\Gamma^{*} : H^{n}_{prim} (X, \qit) \to H^{n-2} (F, \qit(-1))$ is injective. 
\begin{proof} This is likely well-known; however, since a reference for arbitrary $n$ was not found, we proceed as in \cite{BD}. We drop the weights for ease of notation. Then, we consider the projections $p: P \to F$, $q: P \to X$ and note there is a decomposition:
\[ H^{n} (\pit(\E)) = p^{*}H^{n} (F) \oplus  p^{*}H^{n-2} (F)\cdot q^{*}h\] 
From \cite{F} Chapter 3 (or otherwise), we have
\[ p_{*}(p^{*}\gamma_{n}) = 0 , \ p_{*}(p^{*}\gamma_{n-2}\cup q^{*}h) = \gamma_{n-2} \]
for $\gamma_{j} \in H^{j} (F)$. Since $q^{*}$ is injective, it suffices to show that
\[ p^{*}H^{n} (F) \cap  q^{*}H^{n}_{prim} (X) = 0 \]
So, suppose there is some $\gamma \in H^{n}_{prim} (X)$ and some $\gamma_{n} \in H^{n} (F)$ such that
\[ p^{*}\gamma_{n} = q^{*}\gamma \in H^{n} (\pit(\E)) \]
Then we have
\[
p^{*}\gamma_{n-2}\cdot q^{*}h = q^{*} (\gamma \cup h) = 0 \]
since $\gamma\in H^{n}_{prim} (X)$. Pushing forward, it follows that
\[ \gamma_{n-2} = p_{*}(p^{*}\gamma_{n-2}\cup q^{*}h) = p_{*}(q^{*}(\gamma \cup h)) = 0 \]
as was to be shown.
\end{proof}
\end{Lem}
\noindent In particular, this means that there exists some transcendental class 
\[\omega \in \Gamma^{*}H^{n}_{prim} (X, \qit) \subset H^{n-2}(F, \qit(-1)) \]
Thus, (dropping weights) the class
\[ [\prescript{t}{}{\Gamma}] \in H^{2(n-1)} (X \times F, \qit) \]
is not decomposable since its K\"unneth component in
\[ H^{n}_{prim} (X, \qit) \otimes H^{n-2} (F, \qit) \] 
is non-zero. It follows that $\Gamma \not\in \mathcal{D}R^{*}(F \times X)$. To show that this is the case for all non-zero cycles in $\mathcal{D}R^{*}(F \times X)'$, we begin with the following decomposition of the cohomology of $F$. Let
\[ \mathcal{H}:=  H^{n}_{prim} (X, \qit(1))\]
\begin{Thm}[Galkin-Shinder, \cite{GS} Theorem 6.1]\label{GS} There is an isomorphism of Hodge structures:
\[ H^{*} (F, \qit) \cong \text{Sym}^{2} (\mathcal{H}) \oplus \bigoplus_{k=0}^{n-2} \mathcal{H}(-k) \oplus \bigoplus_{k=0}^{2(n-2)} \qit(-k)^{a_{k}} \]
for some positive integers $a_{k}$.
\end{Thm}
\noindent It follows that the transcendental cohomology (not in the middle degree) of $F$ is concentrated in degrees $n-2+2k$ for which $0\leq k \leq n-2$. By the hard Lefschetz theorem, it follows that the K\"unneth component of $\prescript{t}{}{\Gamma}\cdot p^{*}c_{1}^{m}$ in 
\[ H^{n}_{prim} (X, \qit) \otimes H^{n-2+2k} (F, \qit) \]
is non-zero (again, dropping weights). We deduce that
\[ \mathcal{D}R^{*} (F \times X) \cap  \mathcal{D}R^{*} (F \times X)' = 0\]

\end{proof}
\end{Thm}
\begin{Cor}\label{crucial-cor} The cycle class map
\[ R^{k} (F \times X) \to H^{2k} (F \times X, \qit(k)) \]
is injective for all $k$.
\begin{proof} This follows directly from Lemma \ref{Dr-inj} and the proof of Theorem \ref{key}, which shows that the cycle class map
\[ \mathcal{D}R^{*} (F \times X)' \to H^{2*} (F \times X, \qit(*))\]
is injective and that $\mathcal{D}R^{*} (F \times X) \cap \mathcal{D}R^{*} (F \times X)' = 0$ cohomologically.
\end{proof}
\end{Cor}

\section{More Tautological rings}

\noindent We introduce two final tautological rings. Indeed, we let $X^{n}$ denote the $n$-fold product of $X$ and let
\[ \pi_{i}: X^{n} \to X, \ \pi_{ij}: X^{n} \to X^{2}\]
denote the projections onto the $i^{th}$ and $(i,j)^{th}$ factors. For $n=3$ let 
\[ \Delta_{ij} = \pi_{ij}^{*}\Delta_{X} \]
and let $\Delta^{3} = \Delta_{12}\cdot \Delta_{23}$ denote the small diagonal. We also retain the notation of $X^{[2]}$ for the Hilbert scheme. Moreover, let $\epsilon: \widetilde{X^{2}} \to X^{2}$ denote the blow-up of $X^{2}$ along the diagonal and $\pi: \widetilde{X^{2}} \to X^{[2]}$ be the natural projection.
\begin{itemize}
\item $R^{*} (X^{3})$, the subring generated by $\pi_{i}^{*}R^{*}(X)$ and the diagonals $\Delta_{12}$, $\Delta_{23}$, and $\Delta_{13}$ (or, equivalently, the subring generated by $A^{1} (X^{3})$ and the diagonals);
\item $R^{*} (X^{[2]} \times X) = (\pi \times \text{id}_{X})_{*}(\epsilon \times \text{id}_{X})^{*}R^{*} (X^{3})$.
\end{itemize}
\noindent Note that by definition we have
\[ \pi_{ij}^{*}R^{*} (X^{2}) \subset R^{*} (X^{3}) \]
\indent Now, we will consider the correspondence from Theorem \ref{Later} (or rather its transpose). Explicitly, this is given as follows. Let
\[ i: P^{[2]}:= \pit(\Sym^{2}(\E)) \to X^{[2]} \]
be the inclusion of 
\[ Z = \{ x+y \in X^{[2]} \ | \ \exists \ell \in F \text{ such that } x, y \in \ell \} \]
and let $p^{[2]}: P^{[2]} \to F$ be the projection. Then,
\[ \Phi:= \Gamma_{p^{[2]}}\circ\prescript{t}{}{\Gamma_{i}} \in CH^{2n-2} (X^{[2]} \times F) \]
induces (the transpose of) the map $\h(X^{[2]}) \to \h(F)(-2)$ in Theorem \ref{Later}. There is also the fiber diagram:
\[ \begin{tikzcd}P^{2} \arrow{r}{i'} \arrow{d}{\pi'} &\widetilde{X^{2}} \arrow{d}{\pi}\\
P^{[2]} \arrow{r}{i}  & X^{[2]}
\end{tikzcd}\]
where $i': P^{2} = P \times_{F} P \to \widetilde{X^{2}}$ is induced by the universal property of blow-up from the natural map
\[ \rho: P^{2} \hookrightarrow P \times P \xrightarrow{q \times q} X \times X \]
Also, let $p^{2}: P^{2} \to F$ be the projection.
\begin{Lem} $(\Phi \times \text{id}_{X})_{*}(\pi \times \text{id}_{X})_{*}(\epsilon \times \text{id}_{X})^{*} = (p^{2} \times \text{id}_{X})_{*}(\rho\times \text{id}_{X})^{*}$
\begin{proof} This follows from the excess intersection formula \cite{F} \S 6.3. Indeed, by the above fiber diagram, we have
\[ (i\times \text{id}_{X})^{*}(\pi \times \text{id}_{X})_{*} = (\pi' \times \text{id}_{X})_{*}(i' \times \text{id}_{X})^{*} \]
Since $\epsilon\circ i' = \rho$ and $p^{[2]}\circ\pi' = p^{2}$, the result now follows.
\end{proof}
\end{Lem}
\begin{Prop}\label{key-prop} $(\Phi \times \text{id}_{X})_{*}R^{*} (X^{[2]} \times X) \subset R^{*-2} (F \times X)$
\begin{proof} By the previous lemma, it suffices to show that
\[(p^{2} \times \text{id}_{X})_{*}(\rho\times \text{id}_{X})^{*}R^{*} (X^{3}) \subset R^{*} (F \times X) \]
For this we first show that 
\begin{equation} (\rho\times \text{id}_{X})^{*}R^{*} (X^{3}) \subset R^{*} (P^{2} \times X)\label{nice} \end{equation}
We recall from the beginning of the previous section that $R^{*} (P^{2} \times X)$ is the subring of $A^{*} (P^{2} \times X)$ generated by $(p^{2} \times \text{id}_{X})^{*}R^{*} (F \times X)$ and 
\[ h_{1}, h_{2} \in A^{1} (P^{2} \times X) \]
are the two tautological bundles for the $\pit^{1} \times \pit^{1}$-bundle 
\[ p^{2} \times \text{id}_{X}: P^{2} \times X \to F \times X \]
In fact, if (\ref{nice}) holds, it is easy to see from the projective bundle formula that the statement will immediately follow. Now, we can reduce verifying (\ref{nice}) to verifying 
\begin{equation}  (\rho\times \text{id}_{X})^{*}\Delta_{ij} \in R^{*} (P^{2} \times X)\label{nice-red}\end{equation}
since $(\rho\times \text{id}_{X})^{*}$ commutes with $\cdot$ and since 
\[ (\rho\times \text{id}_{X})^{*}A^{1} (X^{3}) \subset A^{1} (P^{2} \times X) = R^{1} (P^{2} \times X)\]
We first verify (\ref{nice-red}) for $i=1$, $j=2$. In this case, there is a commutative diagram:
\[\begin{tikzcd}
P^{2} \times X \arrow{r}{\rho \times \text{id}_{X}} \arrow{d}{} & X^{3} \arrow{d}{\pi_{12}}\\
P^{2} \arrow{r}{\rho} & X^{2}
\end{tikzcd}
\]
so that 
\[ (\rho\times \text{id}_{X})^{*}\Delta_{12} \in R^{*} (P^{2} \times X)  \]
since $R^{*} (P^{2} \times X)$ contains the pull-back of $R^{*} (P^{2})$ by construction (and the projective bundle formula). For $i=1,2$, $j=3$, we consider the composition
\[ \rho_{i}:  P^{2}  \xrightarrow{\rho} X^{2} \xrightarrow{\pi_{i}} X\]
Then, there is a commutative diagram:
\[\begin{tikzcd}
P^{2} \times X \arrow{dr}{}[swap]{\rho_{i} \times \text{id}_{X}} \arrow{r}{\rho \times \text{id}_{X}} & X^{3} \arrow{d}{\pi_{i3}}\\
& X^{2}
\end{tikzcd}
\]
from which it follows that
\begin{equation} (\rho \times \text{id}_{X})^{*}\Delta_{i3} = (\rho_{i} \times \text{id}_{X})^{*}\Delta_{X} = \Gamma_{\rho_{i}} \in A^{n} (P^{2} \times X)\label{imp-obs} \end{equation}
Now, we can view $\Gamma_{\rho_{i}}$ as a divisor on 
\[i_{P'}: P' := (p^{2} \times \text{id}_{X})^{-1}(P) \hookrightarrow P^{2} \times X,\] 
and we can view $P'$ as a $\pit^{1}$-bundle over $P^{2}$ for which $c_{1}(\mathcal{O}_{P}(1)) = (q')^{*} (h)$, where $q': P'\xhookrightarrow{i_{P'}} P^{2} \times X \xrightarrow{\pi_{X}} X$ denotes the projection. We also consider
\[ p': P' \xhookrightarrow{i_{P'}} P^{2} \times X \xrightarrow{\pi_{P^{2}}} P^{2} \] 
the projections onto $P^{2}$ and $F$, resp. Then, we have the lemma below:
\begin{Lem} $\Gamma_{p_{i}} = (p')^{*}\rho_{i}^{*}h + (q')^{*}h \in A^{1} (P')$.
\begin{proof} We have the projective bundle formula:
\[ A^{1} (P') = (p')^{*}A^{1} (P^{2}) \oplus \qit\cdot (q')^{*}h  \]
which means that 
 \[ \Gamma_{p_{i}} =  (p')^{*}\alpha + c\cdot (q')^{*}h \]
 for $c \in \qit$ and $\alpha \in A^{1} (P^{2})$. Explicitly, we have 
 \[ \begin{split}  \alpha= p'_{*}(\Gamma_{\rho_{i}}\cdot (q')^{*}h) &= \pi_{P^{2}*}i_{P'*}(\Gamma_{\rho_{i}}\cdot i_{P'}^{*}\pi_{X}^{*}h)\\ &= \pi_{P^{2}*}(\Gamma_{\rho_{i}}\cdot \pi_{X}^{*}h)\\ & = \rho_{i}^{*}h  \end{split}\]
 Moreover, we see that $\pi_{P^{2}*}\Gamma_{p_{i}} =[P^{2}] \in CH^{0} (P^{2})$, from which we deduce that $c=1$. This gives the lemma.
\end{proof}
\end{Lem}
\noindent Using this lemma, (\ref{imp-obs}) becomes
\[\begin{split}(\rho \times \text{id}_{X})^{*}\Delta_{i3} = i_{P'*}((p')^{*}\rho_{i}^{*}h + (q')^{*}h) &= (\pi_{P^{2}}^{*}\rho_{i}^{*}h+ \pi_{X}^{*}h)\cdot[P']\\ &=(\pi_{P^{2}}^{*}\rho_{i}^{*}h+ \pi_{X}^{*}h)\cdot(p^{2} \times \text{id}_{X})^{*}(\Gamma)\\ &\in R^{n} (P^{2} \times X) \end{split}\]
which is the desired result.
\end{proof}
\end{Prop}
\begin{Cor} The cycle class map
\[ R^{k} (X^{[2]} \times X) \to H^{2k} (X^{[2]} \times X, \qit(k)) \]
is injective for all $k$.
\begin{proof} By Theorem \ref{Later} there is an isomorphism of Chow motives:
\[ \h(X^{[2]}) \xrightarrow{\Phi_{*} \oplus \Psi_{*}} \h(F)(-2) \oplus \bigoplus_{i=0}^{n} \h(X)(-i) \]
for which $(\Psi \times \text{id}_{X})_{*}R^{*} (X^{[2]} \times X) \subset R^{*} (F \times X)$ by Proposition \ref{key-prop}. Now, by Corollary \ref{crucial-cor} and Lemma \ref{basic} the cycle class map restricted to $R^{*} (F \times X)$ and $R^{*} (X^{2})$ is injective. Thus, by Lemma \ref{basic} what remains is to show that
\begin{equation} \Psi_{*}R^{*} (X^{[2]} \times X) \subset \bigoplus_{i=0}^{n} R^{*-i} (X \times X)\label{final-goal} \end{equation}
In the notation from earlier, let $\X_{\eta}$ be the generic fiber of the universal family of cubic hypersurfaces and let $\F_{\eta}$ be the corresponding Fano variety of lines. 
Theorem \ref{Later} applies to give an isomorphism for the Chow motives of the generic fibers:
\[ \h(\X_{\eta}^{[2]}) \xrightarrow{\Phi_{*} \oplus \Psi_{*}} \h(\F_{\eta})(-2) \oplus \bigoplus_{i=0}^{n} \h(\X_{\eta})(-i) \]
By definition, it is easy to see that
\[ R^{*} (X^{3}) \subset \text{Im}\{A^{*} (\X_{\eta}^{3}) \to A^{*} (X^{3}) \} \]
and hence that
\[ R^{*} (X^{[2]} \times X) \subset \text{Im}\{ A^{*} (\X_{\eta}^{[2]} \times_{\eta} \X_{\eta}) \to R^{*} (X^{[2]} \times X) \} \]
Certainly, we have 
\[\Psi_{*}A^{*} (\X_{\eta}^{[2]} \times_{\eta} \X_{\eta}) \subset \bigoplus_{i=0}^{n} A^{*-i} (\X_{\eta} \times_{\eta} \X_{\eta}) \]
By specialization, (\ref{final-goal}) then follows from Lemma \ref{first}. 
\end{proof}
\end{Cor}

\section{Proof of Theorem \ref{int-prod}}
\noindent With the notation from the previous section, consider the following truncated modified diagonal diagonal cycle
\[ \gamma^{3} = \Delta^{3} - \frac{1}{3}(\Delta_{12}\cdot \pi_{3}^{*}h^{n} + \text{perm.}) \in R^{2n} (X^{3}) \]
\begin{Lem}\label{hom-lem} There exist $a_{ijk} \in \qit$ for which
\[ [\gamma^{3}] = \sum_{i+j+k=2n} a_{ijk}\pi_{1}^{*}h^{i}\cdot\pi_{2}^{*}h^{j}\cdot\pi_{3}^{*}h^{k} \in H^{4n} (X^{3}, \qit(n)) \]
i.e., $\gamma^{3}$ is decomposable cohomologically.
\begin{proof} Since $\Gamma^{3} \in R^{n} (X^{3})$ it suffices to verify this for the very general cubic hypersurface. To this end, note that the K\"unneth components of $H^{2n} (X^{3}, \qit)$ are all of the form 
\[ H^{i} (X, \qit) \otimes H^{j} (X, \qit) \otimes H^{k} (X, \qit) \]
for $i+j+k=4n$ (where we drop the weights for ease of notation). Since $X$ is very general, $H^{n}_{prim} (X, \qit)$ is an irreducible Hodge structure. So, such summands contain non-trivial Hodge classes only in the following three instances:
\[ (i,j,k) = (n,n,2n), \ (n,2n, n), \ (2n,n,n) \]
Since the only Hodge class in
\[ H^{n}_{prim} (X) \otimes H^{n}_{prim}(X) \cong \End(H^{n}_{prim}(X)) \]
comes from the diagonal, it follows that the only Hodge classes in $H^{4n} (X^{3}, \qit)$ are of the form
\[ \Delta_{12}\cdot\pi_{3}^{*}h^{n} , \ \Delta_{23}\cdot\pi_{1}^{*}h^{n}, \ \Delta_{13}\cdot\pi_{2}^{*}h^{n}, \ \pi_{1}^{*}h^{i}\pi_{2}^{*}h^{j}\pi_{3}^{*}h^{2n-i-j} \]
from which it follows that $\Delta^{3}$ is a linear combination of these Hodge classes. It is an exercise to the reader to verify that the coefficient of $\Delta_{12}\cdot \pi_{3}^{*}h^{n}$ (and permutations) is $\frac{1}{3}$.
\end{proof}
\end{Lem}
\noindent Now, let $a_{ijk} \in \qit$ be as in Lemma \ref{hom-lem} and consider the corresponding augmented modified diagonal cycle:
\[ \Gamma^{3} = \gamma^{3} - \sum_{i+j+k=2n} a_{ijk}\pi_{1}^{*}h^{i}\cdot\pi_{2}^{*}h^{j}\cdot\pi_{3}^{*}h^{k} \in R^{2n} (X^{3}) \]
\begin{Cor}\label{final-cor} $\Gamma^{3} = 0 \in  R^{2n} (X^{3})$.
\begin{proof} By definition, $\gamma^{3}$ is invariant under the natural action of the symmetric group $\mathfrak{S}_{3}$ on $R^{n} (X^{3})$. Since $\Gamma^{3}$ is cohomologically trivial by assumption, it is a straightforward exercise to check that there is also an $\mathfrak{S}_{3}$ symmetry for the coefficients $a_{ijk}$ so that 
\[ \Gamma^{3} \in  R^{2n} (X^{3})\] 
is also $\mathfrak{S}_{3}$-invariant. Thus, $\Gamma^{3}$ vanishes if and only if 
\[ (\pi \times \text{id}_{X})_{*}(\epsilon \times \text{id}_{X})^{*}(\Gamma^{3}) \in R^{2n} (X^{[2]} \times X) \]
vanishes. But since $\Gamma^{3}$ is cohomologically trivial, so is  $(\pi \times \text{id}_{X})_{*}(\epsilon \times \text{id}_{X})^{*}(\Gamma^{3})$. Finally, by Corollary \ref{crucial-cor}, $(\pi \times \text{id}_{X})_{*}(\epsilon \times \text{id}_{X})^{*}(\Gamma^{3})$ vanishes, which is the desired result.
\end{proof}
\end{Cor}
\noindent To complete the proof, let $\alpha \in A^{i} (X)$, $\beta \in A^{j} (X)$ with $i,j>0$ and since $X$ is rationally connected, we will also assume that $i+j < n$. We would like to show that
\begin{equation} \alpha\cdot\beta \in \qit\cdot h^{i+j} \in A^{i+j} (X)\label{int-final} \end{equation}
To this end, we begin with the following lemma:
\begin{Lem} $\alpha\cdot \beta = \pi_{3*}(\pi_{1}^{*}\alpha\cdot\pi_{2}^{*}\beta\cdot\gamma^{3})$
\begin{proof} We compute
\[ \pi_{3*}(\pi_{1}^{*}\alpha\cdot\pi_{2}^{*}\beta\cdot\Delta^{3}) = \Delta_{X}^{*}(\pi_{1}^{*}\alpha\cdot\pi_{2}^{*}\beta) = \alpha \cdot \beta  \in A^{i+j} (X)\]
On the other hand, by the projection formula
\[ \begin{split} \pi_{3*}(\pi_{1}^{*}\alpha\cdot\pi_{2}^{*}\beta\cdot\Delta_{12}\cdot\pi_{3}^{*}h^{n}) & = \pi_{3*}(\pi_{12}^{*}\Delta_{X*}(\alpha\cdot\beta)\cdot\pi_{3}^{*}h^{n})\\ & = \pi_{3*}(\pi_{12}^{*}\Delta_{X*}(\alpha\cdot\beta))\cdot h^{n} = 0 \in A^{i+j} (X) \end{split}\] 
(note $i+j<n$ by assumption). Finally, one has
\[ \pi_{1}^{*}\alpha\cdot\pi_{2}^{*}\beta\cdot\Delta_{23}\cdot\pi_{1}^{*}h^{n}, \  \pi_{1}^{*}\alpha\cdot\pi_{2}^{*}\beta\cdot\Delta_{13}\cdot\pi_{2}^{*}h^{n} = 0\]
(note $i, j>0$ by assumption). The lemma now follows.
\end{proof}
\end{Lem}
\noindent To obtain (\ref{int-final}), note that by the previous lemma and Corollary \ref{final-cor} we have
\[\begin{split} \alpha\cdot\beta & = \pi_{3*}(\pi_{1}^{*}\alpha\cdot\pi_{2}^{*}\beta\cdot\gamma^{3})\\ 
& = \sum_{r+s+t=2n} a_{rst}\pi_{3*}(\pi_{1}^{*}(h^{r}\cdot\alpha)\cdot\pi_{2}^{*}(h^{s}\cdot\beta)\cdot\pi_{3}^{*}h^{t})\\
& = a\pi_{3*}(\pi_{1}^{*}(h^{n-i}\cdot\alpha)\cdot\pi_{2}^{*}(h^{n-j}\cdot\beta)\cdot\pi_{3}^{*}h^{i+j})\\
& = am_{\alpha}m_{\beta}\cdot h^{i+j}\end{split}\]
where $a = a_{rst} \in \qit$ is the coefficient in the case that $r=n-i$, $s=n-j$ and $t = i+j$, $m_{\alpha}, m_{\beta} \in \qit$ are the degrees of the zero-cycles $\alpha\cdot h^{n-i}$ and $\beta\cdot h^{n-j}$. Note that the penultimate line follows from the fact for all other summands in the preceding line, either $r> n-i$ or $s >n-j$ or $t\neq i+j$.

\end{document}